\documentclass[11pt]{amsart}
\usepackage{amsxtra}
\theoremstyle{plain}
\newcommand{\cleqn}{\setcounter{equation}{0}}
\newcommand{\clth}{\setcounter{theorem}{0}}
\newcommand {\sectionnew}[1]{\section{#1}\cleqn\clth}
\newtheorem{theorem}{Theorem}[section]
\newtheorem{lemma}[theorem]{Lemma}
\newtheorem{definition-theorem}[theorem]{Definition-Theorem}
\newtheorem{proposition}[theorem]{Proposition}
\newtheorem{corollary}[theorem]{Corollary}
\newtheorem{definition}[theorem]{Definition}
\newtheorem{example}[theorem]{Example}
\newtheorem{remark}[theorem]{Remark}
\newtheorem{notation}[theorem]{Notation}
\newcommand \bth[1] { \begin{theorem}\label{t#1} }
\newcommand \ble[1] { \begin{lemma}\label{l#1} }

\newcommand \bpr[1] { \begin{proposition}\label{p#1} }
\newcommand \bco[1] { \begin{corollary}\label{c#1} }
\newcommand \bde[1] { \begin{definition}\label{d#1}\rm }
\newcommand \bex[1] { \begin{example}\label{e#1}\rm }
\newcommand \bre[1] { \begin{remark}\label{r#1}\rm }

\newcommand \bnota[1] { \begin{notation}\label{n#1}\rm }
\newcommand {\eth} { \end{theorem} }
\newcommand {\ele} { \end{lemma} }

\newcommand {\epr} { \end{proposition} }
\newcommand {\eco} { \end{corollary} }
\newcommand {\ede} { \end{definition} }
\newcommand {\eex} { \end{example} }
\newcommand {\ere} { \end{remark} }

\newcommand {\enota} { \end{notation} }
\newcommand \thref[1]{Theorem \ref{t#1}}
\newcommand \leref[1]{Lemma \ref{l#1}}
\newcommand \prref[1]{Proposition \ref{p#1}}
\newcommand \coref[1]{Corollary \ref{c#1}}

\newcommand \exref[1]{Example \ref{e#1}}
\newcommand \reref[1]{Remark \ref{r#1}}
\newcommand \lb[1]{\label{#1}}
\def \top {\mathrm{top}}        %% top degrees/ranks and wedges
\def \topwedge {\wedge^\top}
\def \Rset {{\mathbb R}}         %mathsets

\def \L  {{\mathcal{L}}}
\def \al {\alpha}
\def \be {\beta}

\def \Ga {\Gamma}

\def \la {\langle}
\def \ra {\rangle}
\def \hra {\hookrightarrow}

\def \Id { {\mathrm{Id}} }

\def \Span { {\mathrm{Span}} }
\def \sign { {\mathrm{sign}} }

\def \Tr { {\mathrm{Tr}} }
\def \g  {\mathfrak{g}}   % Lie algebra letters
\def \gl  {\mathfrak{gl}}

\def \f  {\mathfrak{f}}

\def \q  {\mathfrak{q}}
\def \p  {\mathfrak{p}}
\def \sl {\mathfrak{sl}}
\DeclareMathOperator \Mod { {\mathrm{Mod}} }
\DeclareMathOperator \ad { {\mathrm{ad}} }
\DeclareMathOperator \Ker { {\mathrm{Ker}} }

\renewcommand \Im { {\mathrm{Im}} }

\begin{document}
%%%%%%%%%%%%%%%%%%%%%%    Title    %%%%%%%%%%%%%%%%%%%%%%%%%%%%%%%%%%%%%%%%
\title[Modular classes]
{Modular classes of regular twisted Poisson
  structures on Lie algebroids}

\author[Yvette Kosmann-Schwarzbach]{Yvette Kosmann-Schwarzbach}
\address{
Centre de Math\'ematiques Laurent Schwartz \\
\'Ecole Polytechnique \\
F-91128 Palaiseau, France}
\email{yks@math.polytechnique.fr}
\author[Milen Yakimov]{Milen Yakimov}
\address{
Department of Mathematics \\
University of California \\
Santa Barbara, CA 93106, U.S.A.}
\email{yakimov@math.ucsb.edu}
%\thanks{}
\date{}
\begin{abstract} 
We derive a formula for the 
modular class of a Lie algebroid with a regular twisted Poisson structure
in terms of a canonical Lie algebroid 
representation of the image of the
Poisson map. We use this formula to compute the modular classes of 
Lie algebras with a twisted triangular $r$-matrix. 
The special case of $r$-matrices associated to Frobenius 
Lie algebras is also studied.
\end{abstract}
\maketitle
%%%%%%%%%%%%%%%%%%%%   Introduction   %%%%%%%%%%%%%%%%%%%%%%%%%%%%%%%%%%%%%%%%
\sectionnew{Introduction}\lb{intro}
Twisted Poisson structures on manifolds, whose definition we recall 
in Section \ref{recall}, first appeared in the
mathematical physics literature. They 
were introduced in geometry by Klim\v{c}\'ik
and Strobl \cite{KS}, and were studied by \v Severa and Weinstein 
\cite{SW} who proved that they can be described as Dirac structures in
Courant algebroids.
Roytenberg then showed in \cite{R} 
that, more generally, twisted Poisson structures on Lie algebroids
appear in a natural way in his general theory of
twisting of Lie bialgebroids.

While the modular vector fields of Poisson
manifolds were first defined by Koszul in 1985, the theory of
the modular classes of Poisson manifolds was developed by Weinstein
in his 1997 article \cite{W} and that of the modular classes of Lie algebroids,
generalizing those of Poisson manifolds, by Evens,
Lu and Weinstein in 1999 \cite{ELW}.

In \cite{KL}, Kosmann-Schwarzbach and Laurent-Gengoux defined the modular class
of a Lie algebroid $A$ with a twisted Poisson structure $(\pi, \psi)$. 
(The untwisted case where $\psi =0$ corresponds to a triangular Lie
bialgebroid, see \cite{K}.) The modular class, which we shall denote by $\theta(A,
\pi, \psi)$, is a class in the first cohomology group $H^1(A^*)$ of the
dual Lie algebroid $A^*$.
Furthermore,
Kosmann-Schwarzbach and Weinstein \cite{KW}
showed that this class
is, up to a factor of $\frac{1}{2}$, 
the relative modular class for the morphism $ \pi^\sharp \colon A^* \to A$,
\begin{equation}
\label{rel}
2 \, \theta(A, \pi, \psi) = 
\Mod(A^*) -
(\pi^\sharp)^* (\Mod A) \ .
\end{equation}
Here $\Mod A \in H^1(A)$ and $\Mod(A^*)\in H^1(A^*)$ 
denote the modular classes of the Lie algebroids $A$ and $A^*$
in the sense 
of Evens, Lu and Weinstein \cite{ELW},
and $\pi^{\sharp}$ is the map 
$\alpha \in A^* \mapsto i_{\alpha} \pi \in A$. 
As was observed in \cite{KL}, if $\pi^{\sharp} \colon A^* \to A$ is invertible,
the modular class $\theta(A, \pi, \psi)$ vanishes. This can be seen from 
\eqref{rel}, since, in this case, 
$\pi^{\sharp} \colon A^* \to A$ is an isomorphism 
of Lie algebroids, and the two terms in \eqref{rel} cancel out.

Let us 
call a twisted Poisson structure $(\pi, \psi)$ 
on a Lie algebroid $A$ {\it regular} if $\pi^\sharp$ has constant rank, 
{\it{i.e.}}, if the image $B$ of
$\pi^{\sharp}$ is a vector sub-bundle and thus a Lie subalgebroid of $A$. 
For such structures, 
our main result, \thref{main}, expresses 
$\theta(A, \pi, \psi)$ as the image of the characteristic class of 
the Lie algebroid $B$ with representation on the
top exterior power of the kernel of $\pi^{\sharp}$
defined by the adjoint action. Formula \eqref{formula1} in \thref{main}
takes into account the cancellations which occur in \eqref{rel}, even when
$\pi$ is not an isomorphism, and this formula can be used
to compute modular classes. 
Section \ref{regular} contains the necessary definitions, in particular that of
the {\it canonical representation}
of the image of $\pi^{\sharp}$ on its kernel, and 
two lemmas that lead to the proof
of \thref{main} in Section \ref{main_sect}.

In Section \ref{symplectic} we describe 
the regular twisted Poisson structures 
on Lie algebroids by means of  
a linearization of condition \eqref{pipsi}, extending
and unifying two very
different results. On the one hand, 
\v Severa and Weinstein observed in \cite{SW} that 
a twisted Poisson manifold $(M,\pi, \psi)$
admits a foliation by submanifolds equipped with
a non-degenerate $2$-form $\omega$ which is $\psi$-twisted,
{\it{i.e.}}, such that its differential is equal 
to the pull-back of $- \psi$ under the canonical
injection, generalizing the symplectic foliation of Poisson manifolds.
On the other hand, it was shown by Stolin \cite{St} and by Gerstenhaber and
Giaquinto \cite{GG}
that non-degenerate triangular $r$-matrices on Lie algebras are in one-to-one
correspondence with quasi-Frobenius structures (a property that was
used by Hodges and Yakimov to describe the geometry of triangular 
Poisson Lie groups in terms of symplectic reduction \cite{HY}). 
In \thref{linearize} and its corollaries, we extend both these results 
to the case of regular twisted
Poisson structures on Lie algebroids.

In Section \ref{baseapoint} we deal with the important particular case
where the Lie algebroid is a Lie algebra, 
considered as a Lie algebroid over a point. 
As a consequence of \thref{main},
in \prref{Liealg} and \coref{dual} we obtain two simple formulas 
for the modular class defined
by a {\it twisted triangular $r$-matrix} on a Lie algebra. 
The linearization of the twisted classical Yang-Baxter equation,
obtained in \prref{linearize2} and \coref{cor2}, as a consequence of
\thref{linearize}, is very useful 
for explicit constructions, as shown in the examples of Section
\ref{examples}. The one-to-one correspondence between non-degenerate 
triangular $r$-matrices and quasi-Frobenius structures 
is extended to the twisted structures in
\coref{twistedquasifrobenius}, 
and we further 
study the case of triangular $r$-matrices associated to {\it Frobenius
Lie algebras}, in which case  
we obtain a simple characterization of the modular
class (\prref{Frobenius}).
In Section \ref{examples}, we
illustrate the use of our formulas with several classes
of examples. We explicitly compute the modular classes of the 
Gerstenhaber--Giaquinto generalized Jordanian $r$-matrix
structures on $\sl_n(\Rset)$ \cite{GG}, and we construct twisted triangular
$r$-matrix structures with trivial and with nontrivial
modular classes.
%%%%%%%%%%%%%%%%%%%%%%%%%%%%%%%%%%%%%%%%%%%%%%%%%%%%%%%%%%%%%%%%%%%%%%%%%%
\sectionnew{The modular class of a regular twisted Poisson structure}
\label{regular}
%%%%%%%%%%%%%%%%%%%%%%%
\subsection{Twisted Poisson structures}\label{recall} 
Recall that a Lie algebroid with a twisted Poisson structure is 
a triple, $(A, \pi, \psi)$, where $A$ is a Lie algebroid 
over a base manifold $M$, $\pi$ and $\psi$ are sections of
$\wedge^2 A$ and $\wedge^3 A^*$, respectively, such that $\psi$ is 
$d_A$-closed, where 
$d_A \colon \Ga(\wedge^\bullet A^*) \to \Ga(\wedge^{\bullet + 1} A^*)$ 
is the Lie algebroid cohomology operator of $A$, satisfying
\begin{equation}
\label{pipsi}
\frac{1}{2}[\pi,\pi]_A = (\wedge^3 \pi^\sharp) \psi \ .
\end{equation}
See \cite{M} for an exposition of the general theory of Lie
algebroids, see \cite{SW} for the twisted Poisson manifolds and
\cite{R} \cite{KL} for the twisted Poisson structures on Lie algebroids.
If $\rho \colon A \to TM$ is the {\it {anchor}} of $A$, then the dual vector bundle
$A^*$ is a Lie algebroid with anchor $\rho \circ \pi^\sharp$ and Lie bracket 
\begin{equation}
[\al, \be]_{\pi, \psi}= \L^A_{\pi^\sharp \al} \be - \L^A_{\pi^\sharp \be} \al 
- d_A (\pi(\al, \be)) 
+ \psi(\pi^\sharp \al, \pi^\sharp \be, \, . \, ) \ ,
\label{A*brac}
\end{equation}
for $\alpha$ and $\beta \in \Gamma(A^*)$, 
where $\L^A$ is the Lie derivation defined by $\L^A_X= [i_X,d_A]$, for
each $X \in \Gamma (A)$.
The map $\pi^\sharp \colon A^* \to A$ is 
a morphism of Lie algebroids. 
%%%%%%%%%%%%%%%%%%%%%%%%%%
\subsection{A canonical representation in the regular case} 
\label{2.2}
Let $(\pi, \psi)$ be a twisted Poisson 
structure on the Lie algebroid $A$, and assume that it is regular. 
Let $B \subset A$ be the image of $\pi^{\sharp}$.
The morphism of Lie algebroids 
$\pi^\sharp \colon A^* \to A$ 
and the canonical embedding of Lie algebroids
$\iota_B \colon B \hra A$ define
maps of cohomology spaces,
\begin{equation*}
(\pi^\sharp)^* \colon H^\bullet(A) \to H^\bullet(A^*) \quad
{\mathrm{and}} \quad
(\iota_B)^* \colon H^\bullet(A) \to H^\bullet(B) \ .
\end{equation*}
Since the image of $\pi^\sharp \colon A^* \to A$
is $B$, there is an induced map, 
$$
\pi^{\sharp}_B \colon A^* \to B \ ,
$$
which in turn defines a map of cohomology spaces,
\begin{equation*}
\label{piB}
(\pi^\sharp_B)^* \colon H^\bullet(B) \to H^\bullet(A^*) 
\end{equation*}
satisfying
\begin{equation}
\label{comp}
(\pi^\sharp)^* = (\pi^\sharp_B)^* \, \circ \, \iota_B^* \ .
\end{equation}

Let $C \subset A^*$ be the kernel of the morphism $\pi^{\sharp}$.
It follows from the definition \eqref{A*brac} of the 
Lie bracket on $\Ga(A^*)$ that  
$\Ga(C)$ is an abelian ideal of $\Ga(A^*)$, and 
there is an exact sequence of Lie algebroids over the same base,
\begin{equation}
\label{e_seq}
0 \to C \to A^* \to B \to 0 \ .
\end{equation}

Because $\pi$ is skew-symmetric, 
$C$ is the orthogonal of $B$ with respect to 
the fiberwise duality pairing, $\la \, . \, , \, . \, \ra$ between $A$
and $A^*$. Therefore the vector bundles $A/B$ and $C^*$
are isomorphic, and the vector bundles $A^*/C$ and $B^*$
are isomorphic.

The Lie algebra $\Gamma (A)$ acts on $\Ga (A^*)$ by the coadjoint action
(relative to the Lie bracket $[\, . \, , \, . \,]_A$), and the
restriction of this action to $\Ga (B)$ leaves $\Ga (C)$ invariant. 
This is proved using the skew-symmetry and morphism properties of
$\pi^{\sharp}$. However
this is not a representation of $B$ on $C$ unless the anchor of $A$
vanishes. 

On the other hand, from \eqref{e_seq}
we see that the adjoint action of $\Ga(A^*)$ on
$\Ga(C)$ (relative to the Lie bracket $[\, . \, , \, . \,]_{\pi,\psi}$)
factors through the submersion
$\pi^{\sharp}_B : A^* \to~B$, and
therefore induces an action of $\Ga(B)$, which is in fact a representation
of $B$ on $C$, because $\rho \, \circ \,
\pi^\sharp$ vanishes on $C$.  It is an example
of the representation attached to an abelian extension of Lie algebroids
(\cite[Proposition 3.3.20]{M}). 
It follows from formula \eqref{A*brac} 
that this action is defined by
\begin{equation}
\label{B_act_ker1}
X \cdot \gamma = \L^A_X ( \gamma) \ , 
\end{equation}
for all $X \in \Ga(B)$ and $\gamma \in \Ga(C)$.
Explicitly, for all $Y \in \Gamma(A)$,
\begin{equation}
\label{B_act_ker2}
\la X \cdot \gamma , Y \ra = \rho(X) \la \gamma, Y \ra - \la \gamma ,
[X,Y]_A \ra \ .
\end{equation}
We shall refer to this representation of the Lie algebroid $B$ on the
vector bundle $C$ as its {\it canonical representation}.

When restricted to $\Ga (B)$, the adjoint action of $\Ga (A)$ on itself  
induces an action of $\Ga(B)$ on $\Ga(A/B)$, which is defined by
\begin{equation}\label{dualrep}
X \cdot cl Y = {cl} ([X,Y]_A) \ ,
\end{equation}
where $X$ is a 
section of $B$, $Y$ a section of $A$, and ${cl}$ denotes the class
of a section of $A$ modulo the sections of $B$.
The right-hand side is 
well defined because $B$ is a Lie subalgebroid of $A$, and the map
$(X,cl Y) \mapsto X \cdot cl Y$ is 
$C^\infty(M)$-linear in $X$. It is therefore clear that
\eqref{dualrep} defines a representation of $B$ on $A/B$.
It is an example
of a Bott representation as defined by Crainic \cite[Examples 4]{C}. 
It is easy to prove that, under the identification of $C^*$ with
$A/B$, this representation is dual to the canonical representation. 
(We recall the definition of dual representations in Section \ref{char} below.)
\bre{remark} 
When the anchor of $A$ vanishes, {\it e.g.}, 
in the case of Lie algebras, 
the action of $\Ga (B)$ on $\Ga(C)$ which is the restriction 
of the coadjoint action of $\Ga(A)$ on $\Ga(A^*)$ 
defines a representation of $B$ on $C$, and 
equation \eqref{B_act_ker2} shows that it coincides with 
the canonical representation 
defined by \eqref{B_act_ker1}.
\ere
%%%%%%%%%%%%%%%%%%%%%%%%%%%%%%%%
\subsection{Characteristic classes of Lie algebroids}\label{char} 
In \cite{ELW}, Evens, Lu and Weinstein defined 
the characteristic class of a Lie algebroid $E$ with a representation 
on a line bundle $L$, over the same base manifold $M$. 
If $L \to M$ is trivial, the characteristic class 
is the class in $H^1(E)$ of the section $\theta_s \in \Ga (E^*)$
such that, for all $x \in \Ga(E)$, 
\begin{equation}
x \cdot s = \la \theta_s, x \ra s, \quad \
\end{equation}
where $s \in \Ga (L)$ is a nowhere-vanishing section. In fact, the section
$\theta_s$ is a $d_E$-cocycle, and its class, which we shall denote by
$\theta (E, L)$, is
independent of the choice of the section $s$.
If $L$ is not trivial, its characteristic class is defined as one-half that of its square.

The {\it modular class of a Lie algebroid} $E$ over $M$
is the characteristic
class of the canonical representation of $E$ on the line bundle 
$\topwedge E \otimes \topwedge (T^*M)$, defined by
\begin{equation}
x \cdot (\lambda \otimes \mu) = \mathcal L^E_x \lambda \otimes \mu + \lambda
\otimes \mathcal L_{\rho(x)} \mu \ ,
\end{equation} 
for all $x \in \Gamma(E)$, where $\lambda$ and $\mu$ are 
nowhere-vanishing sections of
$\wedge^{\rm {top}} E$ and $\wedge^{\rm {top}} T^*M$, respectively.
We denote the modular class of $E$ by $\Mod E$.

We recall that representations of a Lie algebroid $E$ on vector
bundles $V$ and $V^*$ are {\it dual} if
$$
\la x \cdot s, \sigma \ra = - \la s, x \cdot \sigma \ra + \rho(x) \la s, \sigma
\ra \ ,
$$
for all $x \in \Ga ( E)$, $s \in \Ga (V)$ and  $\sigma \in \Ga (V^*)$.

The proof of the following lemma is simple and will be omitted.
\ble{map} (i) If $f \colon A_1 \to A_2$ is a morphism of Lie algebroids 
over the same base
and $V$ is a representation of $A_2$, then
\[
f^* \theta (A_2, V) = \theta (A_1, V^f) \ ,
\]
where $V^f$ is the representation of $A_1$ induced from $f$.

\noindent (ii)  Dual representations have opposite characteristic classes.
\ele
%%%%%%%%%%%%%%%%%%%%%%%%%%%%%%%%%%%
\subsection{Relations between $\Mod(A^*)$, $\Mod A$
and $\Mod B$} 
The following lemmas contain formulas that will be used in the
proof of our main result, \thref{main}.
\ble{plusminus} Let $A$ be a Lie algebroid with a regular twisted Poisson
structure $(\pi, \psi)$; let $B$ be the image of $\pi^{\sharp}$, and
let $C$ be its kernel.
Denote by $\theta_B$ the 
characteristic class of the Lie algebroid $B$ with the representation 
on the line bundle
$\topwedge (C)$ induced by the canonical representation
\eqref{B_act_ker1}. Then
\begin{equation}
\label{eq1}
\Mod(A^*) = 
(\pi^\sharp_B)^* (\Mod B + \theta_B) \ .
\end{equation}
\ele
\begin{proof}
We shall assume that $\topwedge C$, $\topwedge B$,  and 
$\topwedge A^* \otimes \topwedge T^* M$ are trivial bundles over $M$. If 
this is not the case, the proof below is easily modified using densities.
Let $s_1$ be a nowhere-vanishing section of $\topwedge C$.
Let ${\rm{rk}} B$ denote the rank of $B$, and let 
$s_2 \in \Ga ( \wedge ^{{\rm{rk}} B} A^*)$
be such that $s_1 \wedge s_2$ and $\pi^\sharp s_2$ are 
nowhere-vanishing sections of $\topwedge A^*$ and $\topwedge B$, respectively.
Finally let $\mu \in \Ga (\topwedge T^* M)$
be a volume form on $M$.
We denote by $\xi$, $\eta$ and $\zeta$ representatives of $\Mod A^*$,
$\Mod B$ and $\theta (A^*, \topwedge C)$, corresponding
to the sections $s_1 \wedge s_2 \otimes \mu$, $\pi^\sharp s_2 \otimes
\mu$ and $s_1$ of 
$\topwedge A^* \otimes \topwedge T^*
M$, $\topwedge B \otimes \topwedge T^* M$ and
$\topwedge C$, respectively. By definition, for all $\alpha \in \Ga (A^*)$,
\begin{equation*}
<\xi, \alpha> s_1 \wedge s_2 \otimes \mu = [\alpha, s_1 \wedge s_2]_{\pi, \psi}
\otimes \mu + s_1 \wedge s_2 \otimes 
\mathcal L_{\rho \circ \pi^\sharp (\alpha)} \mu \ ,
\end{equation*}
and 
\begin{equation*}
<\zeta, \alpha> s_1 = [\alpha, s_1]_{\pi, \psi} \ .
\end{equation*}
Since 
\begin{equation*}
[\alpha, s_1 \wedge s_2]_{\pi, \psi} = [\alpha, s_1]_{\pi, \psi}
\wedge s_2 + s_1 \wedge
[\alpha, s_2]_{\pi, \psi} \ ,
\end{equation*}
we obtain
\begin{equation*}
<\xi - \zeta, \alpha> s_2 \otimes \mu = [\alpha, s_2]_{\pi, \psi}
\otimes \mu +
s_2 \otimes \mathcal L_{\rho \circ \pi^\sharp (\alpha)} \mu \ ,
\end{equation*}
and applying $\pi_B^\sharp \otimes \Id$ to both sides of the preceding
equality yields 
\begin{equation*}
<\xi - \zeta, \alpha> \pi_B^\sharp s_2 \otimes \mu = \pi_B^\sharp
[\alpha, s_2]_{\pi, \psi}
\otimes \mu + \pi_B^\sharp s_2 
\otimes \mathcal L_{\rho \circ \pi^\sharp (\alpha)} \mu \ .
\end{equation*}
From the definition of $\eta$, we obtain
\begin{equation*}
<\eta, \pi^\sharp \alpha> \pi^\sharp s_2 \otimes \mu = [\pi^\sharp
\alpha, \pi^\sharp s_2]_A
\otimes \mu + \pi^\sharp s_2 \otimes 
\mathcal L_{\rho \circ \pi^\sharp (\alpha)} \mu \ .
\end{equation*}
Formula \eqref{eq1} follows, using the morphism property of
$\pi^\sharp$ and \leref{map} ({\it i}).
\end{proof}
\ble{plusminus2} Under the assumptions of \leref{plusminus},
\begin{equation}
\label{eq2}
\iota_B^* (\Mod A ) =  \Mod B - \theta_B \ .
\end{equation}
\ele
\begin{proof}
Using \leref{map} ({\it i}) we obtain 
\begin{equation}
\label{part2a}
\iota_B^* (\Mod A ) = \theta (B, 
\topwedge A \otimes \topwedge T^* M) \ .
\end{equation}
As in the proof of \eqref{eq1}, one shows that, 
in terms of the representation of $B$ on $\topwedge (A/B)$ induced by
the representation defined by \eqref{dualrep},
\begin{equation}
\label{part2}
\theta (B, \topwedge A \otimes \topwedge T^* M) =
\Mod B + \theta (B, \topwedge (A/B)) \ .
\end{equation}
Since the representation
\eqref{dualrep}
is dual to the canonical representation 
\eqref{B_act_ker1},
by \leref{map} ({\it {ii}}),
\[
\theta (B, \topwedge (A/B)) = - \,  \theta_B \ .
\]
Combining this result with \eqref{part2a} and \eqref{part2} 
implies \eqref{eq2}.
\end{proof}

\subsection{A formula for the modular class} 
\label{main_sect}
Taking into account \eqref{comp}, we see that \eqref{eq2} 
implies 
\begin{equation}
\label{eq3}
(\pi^\sharp)^*( \Mod A ) = 
(\pi^\sharp_B)^* (\Mod B - \theta_B) \ .
\end{equation}
Formulas \eqref{eq1} and \eqref{eq3} express the fact that 
$\Mod(A^*)$ and $(\pi^\sharp)^* (\Mod A)$
differ from $(\pi^\sharp_B)^* (\Mod B)$
by opposite quantities.  Our main result states that the modular class
$\theta(A,\pi,\psi)$ of $(A,\pi,\psi)$, which is the class of a $1$-cocycle of
the Lie algebroid $A^*$, is the pull-back of the class of a
$1$-cocycle of the image of $\pi^{\sharp}$.
\bth{main} Let $A$ be a Lie algebroid with a regular twisted Poisson structure
$(\pi, \psi)$, and let $B = \Im \, \pi^\sharp$ and $C = \Ker
\pi^\sharp$. Then the modular class
$\theta(A,\pi,\psi)$ of $(A,\pi,\psi)$ satisfies
\begin{equation}
\theta(A,\pi,\psi) = (\pi^\sharp_B)^* (\theta_B) \ , 
\label{formula1}
\end{equation}
where $\theta_B$ is the 
characteristic class of the Lie algebroid $B$ with the representation 
on the line bundle
$\topwedge (C)$ induced by the canonical representation \eqref{B_act_ker1}.
\eth
\begin{proof}  
The formula follows from  
\leref{plusminus} and \eqref{eq3} together with \eqref{rel}.
\end{proof}
In the particular case of a regular Poisson manifold,
the fact that the modular class is in the image of $(\pi^\sharp_B)^*$ was 
proved by Crainic \cite[Corollary 9]{C}.

Formula \eqref{formula1} permits an effective
computation of the modular class without computing terms
which mutually cancel in $\Mod(A^*)$ and $(\pi^\sharp)^* (\Mod A)$.
In Section \ref{baseapoint}, 
we shall illustrate this with various examples in the particular case
of Lie algebras, considered as Lie algebroids over a point.

%%%%%%%%%%%%%%%%%%%%%%%%%%%%%%%%%%%%%%%%%%%%%%%%%%%%%%%%%%%%%%%%%%%%%%%%%%
\sectionnew{Non-degenerate twisted Poisson structures}
\label{symplectic}
%%%%%%%%%%%%%%%%%%%%%%%%%%%%%%%%%%%%%%%%%%%%%%%
Let $(A, \pi, \psi)$ be a Lie algebroid with a regular twisted
Poisson structure. As above, we 
denote by $B$ the Lie subalgebroid of $A$ which is 
the image of the Lie algebroid 
morphism, $\pi^\sharp \colon A^* \to A$. 
Because $\pi$ is skew-symmetric, $\pi^{\sharp}$ defines an isomorphism 
from $B^*$ to $B$, where $B^* = A^*/\Ker \pi^{\sharp}$ is the
dual of $B$ equipped with the Lie bracket inherited from that of
$A^*$. We shall denote this isomorphism by $\pi^{\sharp}_{(B)}$.
Let $\pi_{(B)}$ be the corresponding bivector on $B$. Thus $\pi_{(B)}$
is nothing but $\pi$ thought of as an element of 
$\Ga (\wedge^2 B) \subset \Ga(\wedge^2 A)$. Then $(\pi_{(B)},
\iota_B^*\psi)$ is a non-degenerate twisted
Poisson structure on $B$. (Recall that $\iota_B \colon B
\hookrightarrow A$ denotes the inclusion.)
We denote the inverse of $\pi^{\sharp}_{(B)}$ by
$\omega_B^\flat$, and we define the $2$-form $\omega_B$ on $B$ by
\begin{equation}
\lb{mu2}
i_X (\omega_B) = - \omega_B^\flat(X) \ ,
\end{equation}
for $X \in \Gamma (B)$. With these conventions, $\omega = \omega_B$
and $\pi = \pi_{(B)}$ satisfy
\begin{equation}\label{inverses}
\omega(\pi^\sharp \alpha, \pi^\sharp \beta)= \pi(\alpha, \beta) \ ,
\end{equation}
for all $\alpha$ and $\beta \in \Ga(A^*)$.
We call the fiberwise non-degenerate $2$-cochain 
$\omega_B$ and the non-degenerate bivector 
$\pi_{(B)}$ {\it inverses} of one another.

Set 
$\psi_B = \iota_B^* (\psi)$.
Then $\psi_B \in \Ga( \wedge^3 B^*)$ is 
a $3$-cocycle of $B$ for which one easily shows that \eqref{pipsi}
implies 
\begin{equation}
d_B (\omega_B) = - \, \psi_B \ ,
\end{equation}
where $d_B$ denotes the Lie algebroid cohomology 
operator of $B$.

Conversely, assume that $A$ is a Lie algebroid, $B$ is a Lie
subalgebroid,\break
$\omega \in \Ga(\wedge^2 B^*)$ is a fiberwise 
non-degenerate $2$-cochain of $B$, and 
$\psi \in \Ga(\wedge^3 A^*)$ is a $3$-cocycle whose pull-back to $B$
is $-d_B \omega$. Let $\pi \in \Ga(\wedge^2 B)\subset \Ga(\wedge^2 A)$
be the inverse of $\omega$.
One shows by a direct computation 
that $(\pi, \psi)$ is a regular twisted Poisson 
structure on the Lie algebroid $A$, for which $B$ is the image of 
$\pi^\sharp$. This proves the following result, 
which yields a linearization of  
the defining condition \eqref{pipsi} of twisted Poisson structures for
regular Lie algebroids.
\bth{linearize} There is a one-to-one 
correspondence between regular twisted Poisson structures on a Lie algebroid
$A$ and triples consisting of
\begin{enumerate}
\item a Lie subalgebroid $B$ of $A$,
\item a fiberwise non-degenerate $2$-cochain $\omega \in \Ga(\wedge^2 B^*)$, and
\item a $3$-cocycle $\psi \in \Ga(\wedge^3 A^*)$
whose pull-back to $B$ is $-d_B \omega$, 
\end{enumerate}
The twisted Poisson bivector $\pi \in \Ga (\wedge^2 B)$ and the
non-degenerate $2$-cochain $\omega$ of $B$ are inverses of one another.
\eth

\bco{cor1} Assume that $A$ is a Lie algebroid and
$\omega \in \Gamma(\wedge^2 A^*)$ is a $2$-cochain of $A$ whose restriction
$\omega|_B$ to a Lie subalgebroid $B$ of $A$ is fiberwise non-degenerate. If
$\pi \in \Ga(\wedge^2 B)$ denotes the inverse of
$\omega|_B$, then $(\pi , - d_A \omega)$ is a 
regular twisted Poisson structure on the Lie algebroid $A$.
\eco

\coref{cor1} shows how one can construct many regular
twisted Poisson structures on Lie algebroids.

%%%%%%%%%%%%%%%%%%%%%%%%%%%%%%%%%%%%%%%%%%%%%%%%%%%%%%%%%
\sectionnew{Lie algebras with twisted triangular $r$-matrices}
\label{baseapoint}
%%%%%%%%%%%%%%%%%%%%%%%%%%%%
\subsection{Twisted triangular $r$-matrices on Lie algebras}
\label{tt1}
The case of a twisted Poisson structure for a Lie algebroid
over a point yields the concept of a twisted triangular 
$r$-matrix for a Lie algebra. Given a Lie algebra $\g$, we say 
that $(r, \psi)$ is a {\it {twisted triangular $r$-matrix 
structure}} on $\g$, or simply a {\it{twisted triangular structure}}, 
if $\pi = r \in \wedge^2 \g$ and $\psi$ is a $3$-cocycle of $\g$, 
where this data satisfies
\eqref{pipsi}, which is equivalent to
\[
[r_{12}, r_{13}] + [r_{12}, r_{23}] + [r_{13}, r_{23}] = 
- (\wedge^3 r^\sharp) \psi \ ,
\]
the {\it twisted classical Yang-Baxter equation}. Here
as above, $r^\sharp \colon \g^* \to \g$ is the linear map defined by 
$r^\sharp(\alpha) = i_{\alpha} r$. 
A twisted triangular structure $(r, \psi)$  on $\g$ gives rise to the 
Lie algebra structure on $\g^*$ defined by \eqref{A*brac}, which, in
this case, reduces to
\[
[\alpha, \beta]_{r, \psi} = \ad^*_{r^\sharp \alpha}(\beta) - 
\ad^*_{r^\sharp \beta}(\alpha) +
\psi(r^\sharp \alpha, r^\sharp \beta, \, . \, ), 
\quad \alpha, \beta \in \g^* \ .
\]

The map $r^\sharp$
is a homomorphism of Lie algebras. 
Its image, denoted by $\p$, will
be called the {\it carrier} of the twisted triangular structure
$(r, \psi)$, in accordance with the terminology for 
triangular $r$-matrices \cite{GG}. The kernel of $r^\sharp$
is an abelian ideal of $\g^*$, and it is stable under 
the coadjoint action of $\p$. As in the general case of Lie
algebroids, we denote  by $r^{\sharp}_{(\p)}$ 
the isomorphism from  $\p^*$ to $\p$ 
defined by $r$.

For a finite-dimensional representation $V$ of a Lie algebra
$\f$, we will denote by $\chi_{\f, V} \in \f^*$ its 
{\it infinitesimal character}, defined by 
\[
\chi_{\f, V}(x) = \Tr_V(x), \quad x \in \f \ .
\] 
The infinitesimal character of the representation $V$ 
is a $1$-cocycle of
$\f$ whose cohomology class 
is the characteristic class of $\f$ with representation on the
$1$-dimensional vector space, $\topwedge V$.

\subsection{Modular classes of twisted triangular structures}
The modular class of a twisted triangular 
structure $(r, \psi)$ on a Lie algebra $\g$ is the class of
a $1$-cocycle of the
Lie algebra $\g^*$, which can be computed by means of the formulas
in the following proposition and its corollary.

\bpr{Liealg} 
Let $\theta(\g,r, \psi)$
be the modular class of the Lie algebra $\g$ with the
twisted triangular structure $(r, \psi)$, with carrier $\p$. Then
\begin{equation}\label{liealgebra}
\theta(\g,r, \psi) 
= - [r^\sharp_{(\p)} (\chi_{\p, \, \Ker r^\sharp})] \ ,
\end{equation}
where $\chi_{\p, \, \Ker r^\sharp}$
is the infinitesimal character of the coadjoint representation of $\p$
on $\Ker  r^\sharp$, and $[x]$ denotes the class in the Lie algebra
cohomology of $\g^*$ of a $1$-cocycle $x \in \p \subset \g$.
\epr
\begin{proof} 
We have observed in \reref{remark} that
the canonical representation of $\p$ on $\Ker r^{\sharp}$, which was 
defined by \eqref{B_act_ker1} in the general case of
regular twisted Poisson structures on Lie algebroids,  
reduces to the coadjoint representation in the case of Lie algebras. 
Therefore
\prref{Liealg} follows from \thref{main}, after  
observing that $r_{(\p)}^\sharp$ is a skew-symmetric map from $\p^*$
to $\p \subset \g$.
\end{proof}

The Lie algebra $\p$ also acts on $\g/\p$ by the action induced from
the adjoint action of $\g$. This representation is dual to the
coadjoint representation of $\p$ on $\Ker r^\sharp$. Therefore

\bco{dual}
Let $\chi_{\p, \, \g/ \p}$ be the infinitesimal character of the
representation of $\p$ on $\g/\p$ induced from the adjoint action.
Then
\begin{equation}\label{liealgebra2}
\theta(\g,r, \psi) =
[ r^\sharp_{(\p)} (\chi_{\p, \, \g/ \p})] \ .
\end{equation}
\eco 

\subsection{Non-degenerate twisted triangular structures}
\thref{linearize} and \coref{cor1}, 
when applied to
the case of Lie algebroids with base a point, yield the following statements.

\bpr{linearize2} There is a one-to-one correspondence between
twisted triangular $r$-matrix structures on a Lie algebra $\g$ 
and triples consisting of
\begin{enumerate}
\item a Lie subalgebra $\p$ of $\g$,
\item a non-degenerate $2$-cochain $\mu \in \Ga(\wedge^2 \p^*)$, and
\item a $3$-cocycle $\psi \in \Ga(\wedge^3 \g^*)$ whose restriction to
  $\p$ is $-d_{\p} \mu$,
where $d_{\p}$ denotes the Lie algebra
cohomology operator of $\p$.
\end{enumerate}
The twisted $r$-matrix $r \in \wedge^2 \p$ 
and the non-degenerate $2$-cochain $\mu$ of $\p$ are inverses of one
another.
\epr

\bco{cor2} Assume that $\g$ is a Lie algebra and
$\mu \in \Gamma(\wedge^2 \g^*)$ is a $2$-cochain of $\g$ whose
restriction $\mu|_\p$ to a Lie subalgebra $\p$ of $\g$ is non-degenerate.
If $r \in \Ga(\wedge^2 \p)$ denotes the inverse of
$\mu|_\p$, then $(r, - d_\g \mu)$ is a twisted triangular r-matrix
structure on the Lie algebra $\g$.
\eco

\prref{linearize2} shows how the twisted classical Yang-Baxter equation
can be linearized, and 
\coref{cor2} shows how one can construct many
twisted triangular
$r$-matrices. In particular, it explains 
the origin of \exref{twist_triang2} below,
as well as the construction of Example 5 in \cite{KL},
see \exref{twist_triang1}. 

\smallskip

When $\mu$ is a non-degenerate $2$-cocycle of a Lie algebra $\p$ 
with values in the trivial 
representation, the pair $(\p, \mu)$ is said to be a 
{\it quasi-Frobenius Lie algebra}, see \cite{St, GG, HY}. 
Let us say that $(\p, \mu, \psi)$ is a {\it twisted 
quasi-Frobenius Lie algebra} if $\mu$ is a non-degenerate $2$-cochain
of the Lie algebra $\p$ and $\psi$ a $3$-cocycle of $\p$ such that 
$d_\p \mu = - \psi$. We obtain, as a particular case of \prref{linearize2},  
the following correspondence.
\bco{twistedquasifrobenius}
There is a
one-to-one correspondence between  non-degenerate twisted 
triangular structures and twisted quasi-Frobenius structures on
Lie algebras. In this correspondence the non-degenerate twisted
triangular $r$-matrix and
the non-degenerate $2$-cochain are inverses of one another.
\eco
The well-known one-to-one correspondence between non-degenerate 
triangular $r$-matrices and quasi-Frobenius structures on
Lie algebras \cite{St,GG} appears as the particular case 
where $\psi =0$.

\subsection{Frobenius Lie algebras}
\label{Frob}
A quasi-Frobenius Lie algebra structure on $\p$ is defined by a
non-degenerate $2$-cocycle, $\mu$. The special case 
of Frobenius Lie algebras corresponds to the case where $\mu$ is 
the coboundary 
in the Lie algebra
cohomology of $\p$ of some $- \xi \, \in \p^*$.
In other words $(\p, \xi)$, with $\xi \in \p^*$, is a 
{\it Frobenius Lie algebra} 
if the skew-symmetric bilinear form $\mu$, defined by 
\begin{equation}
\label{mu}
\mu(X,Y)= \xi([X,Y]) \ ,
\end{equation}
for $X, Y \in \p$, 
is non-degenerate. 
The bilinear form defined by \eqref{mu} is non-degenerate if and only
if the linear map 
\begin{equation}
\label{isom}
X \in \p \mapsto \ad^*_X(\xi) \in \p^*
\end{equation}
is an isomorphism from $\p$ to $\p^*$.
This implies 
that $(\p, \xi)$ is a Frobenius 
Lie algebra if and only if the coadjoint orbit of $\xi$, under the
action of the
adjoint group of $\p$, is dense in $\p ^*$ (see \cite{HY}).

We now consider the case where $(r, \psi)$ is a twisted triangular
structure on a Lie algebra $\g$ derived from 
a Frobenius structure on the carrier subalgebra $\p$.
The map $r^\sharp \colon \g^* \to \g$ defines a 
linear isomorphism $r^{\sharp}_{(\p)} : \p^* \to \p$. 
Let $\mu \in \wedge^2 \p^*$ be the $2$-cochain on $\p$ inverse of
$r_{(\p)} \in \wedge ^2 \p$,
and assume that $\mu$ is the coboundary of a
$1$-cochain
$- \xi \in \p^*$, $\mu = - d_\p \xi$.
In the following proposition we show how to compute the modular class
of a Lie algebra $\g$ with such a triangular 
$r$-matrix. The notation is that of \prref{Liealg} and \coref{dual}.

\bpr{Frobenius} Let $r$ be a triangular $r$-matrix
on a Lie algebra $\g$ 
derived from a Frobenius structure $\xi$ on the carrier subalgebra $\p$
of $r$. Then the modular class of $(\g,r,0)$ is 
the class in the Lie algebra cohomology of $\g^*$ 
of the unique $1$-cocycle $X \in \p$ such that 
\begin{equation}\label{coadj}
\ad^*_X (\xi) = \chi_{\p, \g/\p} \ .
\end{equation}
\epr
\begin{proof} 
For any $\alpha \in \p^*$, the element  
$Z = r^{\sharp}_{(\p)} (\alpha) \in \p$ satisfies
$\xi( [Z, Y]) = - \alpha (Y)$, for all $Y \in \p$,
and therefore
\[
\ad^*_Z (\xi) = \alpha \ .
\]
Thus $X = r^{\sharp}_{(\p)}(\chi_{\p, \g/\p})$ satisfies
\[
\ad^*_X (\xi) = \chi_{\p, \g/\p} \ .
\]
The uniqueness of $X$ satisfying \eqref{coadj}
follows from the fact that \eqref{isom} is a linear isomorphism. 
By \coref{dual}, the class of $X$ is the modular class of $(\g,r,0)$. 
\end{proof}

\bre{remarkfrobenius}
There are no twisted analogues of the Frobenius Lie algebras because when
$\mu$ is a coboundary, the coboundary of $\mu$ vanishes.
When we apply \thref{linearize} successively to tangent bundles and to Lie
algebras we see that 
 
\noindent 
$\bullet$ 
twisted quasi-Frobenius structures on Lie algebras 
correspond to twisted symplectic structures on
manifolds (non-degenerate $2$-forms),

\noindent 
$\bullet$ 
quasi-Frobenius structures on Lie algebras 
correspond to symplectic structures on
manifolds (non-degenerate closed $2$-forms),

\noindent 
$\bullet$ 
Frobenius structures on Lie algebras 
correspond to exact symplectic structures on
manifolds (non-degenerate exact $2$-forms).
\ere

%%%%%%%%%%%%%%%%%%%%%%%%%%%%%%%%%%%%%%%%%%%%%%%%%%%%%%%%%
\subsection{Examples}\label{examples}
Denote by $e_{ij}$, $1 \leq i,j \leq n$, 
the standard elementary matrices
constituting a basis of $\gl_n(\Rset)$, and by 
$e_{ij}^*$, $1 \leq i,j \leq n$, the dual
basis.
First we re-examine \cite[Example 5]{KL}. The published version of this
article contained an error which we now correct.
\bex{twist_triang1} Denote by $\g$ the Lie subalgebra
of $\gl_3(\Rset)$ spanned by 
$\{e_{ij} \mid 1\leq i \leq 2, 1 \leq j \leq 3\}$.
It is a codimension-one subalgebra of a maximal 
parabolic subalgebra of $\gl_3(\Rset)$, and it
is isomorphic to the Lie algebra of the Lie group 
of affine transformations of $\Rset^2$. The pair $(r,\psi)$, where 
\[
r = e_{11} \wedge e_{22} + e_{13} \wedge e_{23} \quad \quad {\mathrm
  {and}} \quad \quad \psi = - (e_{11}^* + e_{22}^*) \wedge e_{13}^* \wedge
e_{23}^* \ ,
\]
defines a twisted triangular $r$-matrix structure on $\g$,
as shown in \cite[Example 5]{KL}. 
The carrier of $r$ is the subalgebra 
$
\p = 
\Span 
\{ e_{11}, e_{22}, e_{13}, e_{23} \}$.
Then $ \dim (\g / \p) = 2$, and one can easily show 
that $\chi_{\p, \, \g / \p} =0$. Therefore 
$\theta(\g, r, \psi) = 0$. 

We observe that the twisted $r$-matrix 
discussed in this example can be easily 
constructed using \prref{linearize2}.
Define the $2$-cochain on $\g$,
\[
\mu = e_{11}^* \wedge e_{22}^* + e_{13}^* \wedge e_{23}^* \ ,
\]
and compute the coboundary of $- \mu$,
\[
\psi_1 = - d_\g \mu = 
- (e_{11}^* + e_{22}^*) \wedge e_{13}^* \wedge e_{23}^* - 
e_{12}^* \wedge e_{21}^* \wedge e_{22}^* + e_{11}^* \wedge e_{21}^* 
\wedge e_{12}^* \ . 
\]
Then the restriction of $\mu$ to $\p$ is non-degenerate, and 
$r$ is the inverse of $\mu_{|\p}$. \coref{cor2} 
implies that $(r, \psi_1)$ 
is a twisted triangular $r$-matrix structure on $\g$. That the same is true 
for the pair $(r, \psi)$ above 
follows from the fact that $\psi- \psi_1$ is a $3$-cocycle for $\g$,
together with the
equality $(\wedge^3 r^\sharp)(\psi - \psi_1)=0$.
\eex

Next we provide an example of a twisted triangular $r$-matrix
with a nontrivial modular class.
\bex{twist_triang2} Let $\q_{n-1}$ be the Lie subalgebra of $\gl_n(\Rset)$ 
spanned by $\{ e_{ij} \mid 1 \leq i \leq n-1, 1 \leq j \leq n \}$.
It is a codimension-one subalgebra of a maximal parabolic 
subalgebra of $\gl_n(\Rset)$ obtained by omitting the  
$(n-1)$-st negative root of $\gl_n(\Rset)$. 
Define the $2$-cochain on $\gl_n(\Rset)$,
\[
\mu = \sum_{1 \leq i < j \leq n-1} e_{ij}^* \wedge e_{ji}^*
+ \sum_{i=1}^{n-1} e_{ii}^* \wedge e_{in}^* \ .
\]
The $ d_{\gl_n(\Rset)}$-coboundary of $- \mu$ is the following 
$3$-cocycle on $\gl_n(\Rset)$:
\begin{equation*}
\psi =
\sum_{i, j, k=1}^n \sign \, (i-j) 
e_{ik}^* \wedge e_{kj}^* \wedge e_{ji}^*  
+ \sum_{1 \leq i, k \leq n-1,  
i \neq k}
e_{ik}^* \wedge e_{ki}^* \wedge e_{in}^* 
- \sum_{i,k=1}^{n-1} e_{ii}^* \wedge e_{ik}^* \wedge e_{kn}^*
\ ,
\end{equation*}
where, for an integer $n$, $\sign \, n = 1$ if $n > 0$, $\sign \, n = - 1$
if $n < 0$ and $\sign \, 0 =0$.  
The restriction of $\mu$ to $\q_{n-1}$ is non-degenerate,
and we consider its inverse,
\[ 
r =
\sum_{1 \leq i < j \leq n-1} e_{ij} \wedge e_{ji}
+ \sum_{i=1}^{n-1} e_{ii} \wedge e_{in} \in \wedge^2 \q_{n-1}
\subset \wedge^2 (\gl_n(\Rset)) \ .
\]
\coref{cor2} implies that 
$(r, \psi)$ is a twisted triangular $r$-matrix 
structure on $\gl_n(\Rset)$. The carrier of $r$ is 
$\q_{n-1}$. It is easy to compute 
\[
\chi_{\q_{n-1}, \gl_n(\Rset)/\q_{n-1}} =
- \sum_{i=1}^{n-1} e_{ii}^* \ .
\]
Since $r^\sharp e^*_{ii} = e_{in}$, we obtain, from \coref{dual},
\[
\theta(\gl_n(\Rset), r, \psi) = - 
\left[
\sum_{i=1}^{n-1} e_{in}
\right] \ .
\]
\eex

%%%%%%%%%%%%%%%%%%%%%%%%%%%%

\bex{Jordanian} Let $\p_1$ be the 
maximal parabolic subalgebra of $\sl_n(\Rset)$
containing all the upper triangular matrices,
obtained by omitting the first negative root.
Then $(\p_1, \xi)$ is a Frobenius Lie algebra with
\[
\xi = \sum_{i=1}^{n-1} e_{i, i+1}^* \ .
\]
The corresponding triangular $r$-matrix on $\sl_n(\Rset)$
is the Gerstenhaber--Giaquinto generalized Jordanian
$r$-matrix,
\[
r_{GG}= \sum_{k=1}^{n-1} d_k \wedge e_{k,k+1} +
\sum_{i<j} \sum_{m=1}^{j-i -1} e_{i, j-m+1} \wedge e_{j, i+m} \ ,
\] 
where 
\[
d_k= \frac{n-k}{n} ( e_{11} + e_{22} + \ldots + e_{kk}) -
\frac{k}{n} ( e_{k+1, k+1} + e_{k+2,k+2} + \ldots + e_{nn}) \ .
\] 
See \cite{GG}, where Gerstenhaber and Giaquinto 
proved that this
skew-symmetric solution of the classical Yang-Baxter equation is a
boundary point of a subset of the set of solutions of the modified
Yang-Baxter equation containing the semi-classical limit of the
generalization of the quantum $R$-matrices
introduced by Cremmer and Gervais in \cite{CG}, and called it the 
Cremmer--Gervais $r$-matrix.
It is easy to see that 
\[
\chi_{\p_1, \sl_n(\Rset) / \p_1} = - (n-1) e_{11}^* 
+  e_{22}^* + e_{33}^* + \ldots + e_{nn}^* \ .
\]
The unique $X \in \p_1$ that satisfies 
$\ad_X^*(\xi) =  \chi_{\p_1, \sl_n(\Rset) / \p_1}$
is the matrix $X =  - \sum_{k=1}^{n-1} (n-k) e_{k, k+1}$. Therefore, by
\prref{Frobenius},
the modular class of $\sl_n(\Rset)$ with the Gerstenhaber--Giaquinto 
generalized Jordanian $r$-matrix is
\[
\theta(\sl_n(\Rset), r_{GG}, 0) = 
- \left[ \sum_{k=1}^{n-1} (n-k) e_{k, k+1} \right]  \ .
\]
This result constitutes a generalization to all $n$'s of 
\cite[Example 4.2]{KL} which is the case of $n=2$.
\eex

{\bf Acknowledgements.} The second author would like to thank 
\'Ecole Polytechnique for the warm hospitality during the week in
August 2006 
when the preliminary results of this paper were obtained.
The research of M.Y. was partially supported 
by NSF grant DMS-0406057 and an Alfred P. Sloan research fellowship.

%%%%%%%%%%%%%%%%%%%%%%%%%%%%%%%%%%%%%%%%%%%%%%%%%%%%%%%%%%%%%%%%%%%%%%%%%%

%%%%%%%%%%%%%%%%%%%%%%%%%

\begin{thebibliography}{AFMO}
%%%%%%%%%%%%%%%%%%%%
       \bibitem{C} Crainic, M., {\em{Differentiable and algebroid cohomology, 
Van Est isomorphisms, and characteristic classes}},  
Comment. Math. Helv. {\bf{78}}, 681--721 (2003).
       \bibitem{CG} Cremmer, E., and Gervais, J-L., {\em{The quantum
group structure associated with non-linearly extended Virasoro algebras}},
Comm. Math. Phys. {\bf{134}}, 619--632 (1990).  
       \bibitem{E} Elashvili, A. G., {\em{Frobenius Lie algebras}}, 
Funct. Anal. Appl. {\bf{16}}, 94--95 (1982); Engl. transl., 326--328.
       \bibitem{ELW} Evens, S., Lu, J.-H., and Weinstein, A., {\em{Transverse measures,
the modular class and a cohomology pairing for Lie algebroids}},
Quart. J. Math. Ser. 2 {\bf{50}}, 417--436 (1999). 
       \bibitem{GG} Gerstenhaber, M., and Giaquinto, A., {\em{Boundary
solutions of the classical Yang--Baxter equation,}} Lett. Math. Phys.
{\bf{40}}, 337--353 (1997).
       \bibitem{HY} Hodges, T. J., and Yakimov, M., {\em{Triangular Poisson 
structures on Lie groups and symplectic reduction}}, Noncommutative geometry 
and representation theory in mathematical physics, Contemp. Math., 391, 
Amer. Math. Soc., 123--134 (2005).
       \bibitem{KS} Klim\v{c}\'ik, C., and Strobl, T., 
{\em{WZW-Poisson manifolds}}, J. Geom. Phys. {\bf{43}}, 341--344 (2002).
       \bibitem{K} Kosmann-Schwarzbach, Y., {\em{Modular vector fields and
       BV-algebras}}, Poisson Geometry, J. Grabowski and
       P. Urbanski, eds., Banach Center Publications
       {\bf{51}}, 109--129 (2000). 
       \bibitem{KL} Kosmann-Schwarzbach, Y., and Laurent-Gengoux, C.,
{\em{The modular class of a twisted Poisson structure}}, 
Travaux math\'ematiques (Luxembourg) {\bf {16}}, 315--339 (2005).
       \bibitem{KW} Kosmann-Schwarzbach, Y., and Weinstein, A., 
{\em{Relative modular classes of Lie algebroids}}, C. R. Acad. Sci. Paris,
Ser. I {\bf{341}}, 509-514 (2005).
       \bibitem{M} Mackenzie, K. C. H., {\em{General theory of Lie
       groupoids and Lie algebroids}}, London Mathematical Society
       Lecture Note Series 213, Cambridge University Press (2005).
       \bibitem{R} Roytenberg, D., {\em{Quasi-Lie bialgebroids and twisted
       Poisson manifolds}}, Lett. Math. Phys. {\bf{61}},
       123--137 (2002).
       \bibitem{SW} \v Severa, P., and Weinstein, A., {\em{Poisson geometry 
with a $3$-form background}}, Noncommutative geometry 
and string theory (Yokohama, 2001), Progr. Theoret. Phys. 
Suppl. {\bf{144}}, 145--154 (2001).
       \bibitem{St} Stolin, A., {\em{On rational solutions of Yang--Baxter 
equation for $\sl(n)$,}} Math. Scand. {\bf{69}}, 57--80 (1991).
       \bibitem{W} Weinstein, A., {\em{The modular automorphism group
       of a Poisson manifold}}, Comm. Math. Phys. {\bf{23}}, 379--394 (1997).
%%%%%%%%%%%%%%%%%%%%%%%%%%%%%%%%%%%%%%%%%%%%%%%%%%%%%%
\end{thebibliography}
\end{document}